\documentclass[12pt]{article}
\usepackage{amsmath,amsthm,amssymb,latexsym,hyperref,url,graphicx,amsthm,multirow,color,tikz,enumitem,appendix}
\usetikzlibrary{trees}

\newtheorem{theorem}{Theorem}
\newtheorem{conj}[theorem]{Conjecture}
\newtheorem{lemma}[theorem]{Lemma}

\newtheorem{innercustomthm}{Theorem}
\newenvironment{customthm}[1]
  {\renewcommand\theinnercustomthm{#1}\innercustomthm}
  {\endinnercustomthm}

\newcommand{\red}{\mathrm{red}}

\begin{document}

\title{Ascent sequences and the binomial convolution of Catalan numbers}
\author{Lara K. Pudwell\\
\small Department of Mathematics and Statistics \\[-0.8ex]
\small Valparaiso University \\[-0.8ex]
\small Valparaiso, IN 46383 \\[-0.8ex]
\small \texttt{Lara.Pudwell@valpo.edu}
}

\maketitle

\begin{abstract}
In this paper, we consider two sets of pattern-avoiding ascent sequences: those avoiding both 201 and 210 and those avoiding 0021.  In each case we show that the number of such ascent sequences is given by the binomial convolution of the Catalan numbers.  The result for $\{201, 210\}$-avoiders completes a family of results given by Baxter and the current author in a previous paper.  The result for 0021-avoiders, together with previous work of Duncan, Steingr\'{i}msson, Mansour, and Shattuck, completes the Wilf classification of single patterns of length 4 for ascent sequences.
\end{abstract}

\section{Introduction}\label{Introduction}

Given an integer string $x_1\dotsm x_n$, an \emph{ascent} is position $j$ such that $x_j<x_{j+1}$.  Write $\mathrm{asc}(x_1 \dotsm x_n)$ for the number of ascents in $x_1\dotsm x_n$.  An \emph{ascent sequence} $x_1\dotsm x_n$ is a sequence of non-negative integers such that
\begin{enumerate}
\item $x_1=0$, and
\item for $1<i \leq n$, $x_i \leq \mathrm{asc}(x_1\dotsm x_{i-1})+1$.
\end{enumerate}
For example, 01234, 0120102, and 01013 are all ascent sequences, while 01024 is not since $\mathrm{asc}(0102)=2$.  Ascent sequences have been an increasingly frequent topic of study since Bousquet-M\'{e}lou, Claesson, Dukes, and Kitaev related them to $(2+2)$-free posets and enumerated the total number of ascent sequences \cite{BDCK10}, thus also enumerating $(2+2)$-free posets (equivalently, interval orders).  Since then, various authors have connected ascent sequences to a number of other combinatorial objects \cite{CL11,DRKS11,DP10,KR11}; also see \cite[Section 3.2.2]{K11} for additional references.  The number of ascent sequences of length $n$ is given by the Fishburn numbers, Online Encyclopedia of Integer Sequences (OEIS) sequence A022493.

Given a string of integers $x=x_1\dotsm x_n$, the \emph{reduction} of $x$, denoted $\red(x)$ is the string obtained by replacing the $i$th smallest digits of $x$ with $i-1$.  For example, $\mathrm{red}(273772) = 021220$.  A \emph{pattern} is merely a reduced string.  We say that $x=x_1 \dotsm x_n$ \emph{contains} a pattern $p=p_1\dotsm p_k$ if there exists a subsequence of $x$ order-isomorphic to $p$, i.e., there exist indices $1 \leq i_1 < i_2 < \dotsm < i_k \leq n$ such that $\red(x_{i_1}x_{i_2}\dotsm x_{i_k}) = p$.  This is analogous to the classical definition of patterns for permutations, but here patterns may contain repeated digits, and patterns are normalized so that their smallest digit is 0 rather than 1. We write $\mathcal{A}(n)$ for the set of ascent sequences of length $n$ and $\mathcal{A}_B(n)$ for the set of ascent sequences of length $n$ avoiding all patterns in list $B$.  Also, we let $\mathrm{a}_B(n) = \left|\mathcal{A}_B(n)\right|$.

Pattern avoidance in ascent sequences was first studied by Duncan and Steingr\'{i}msson \cite{DS11}.  They focused on avoiding a single pattern of length at most 4 and conjectured relationships between sequences avoiding 201, 210, 0123, 0021, or 1012 and other entries in the OEIS \cite{OEIS}. 

Mansour and Shattuck \cite{MS14} later computed the number of sequences avoiding 1012 or 0123 and showed that certain statistics on 0012-avoiding ascent sequences are equidistributed with other statistics on the set of 132-avoiding permutations.  Callan, Mansour, and Shattuck also identified the complete equivalence class of pairs of length-4 patterns such that $\mathrm{a}_{\sigma,\tau}(n)$ is given by the Catalan numbers in \cite{CMS14}.

In \cite{BP14} Baxter and the present author considered the enumeration of ascent sequences avoiding a pair of patterns of length 3.  In particular there are at least 35 different sequences that can be obtained by avoiding a pair of patterns of length 3 in the ascent sequences context, 16 of which are already known in the OEIS for other combinatorial reasons.  One of these results is the following.

\begin{theorem}[\cite{BP14}, Proposition 16]
\label{BPthm1}
$\left|\mathcal{A}_{201,210}(n)\right|=\sum_{k=0}^{n-1}\binom{n-1}{k}C_k$ for $n \geq 1$.
\end{theorem}

In \cite{BP14} this result is followed by the comment ``We defer the proof itself for a separate paper, however, as it is signficantly more complicated than the arguments above.''  In Section \ref{pair} we give the proof of Theorem 1.    The proof begins with a generating tree which is used to derive a system of functional equations for a family of multivariate generating functions.  The solution to the system is difficult to determine directly, but we experimentally conjecture the solution and then validate that it is indeed the unique solution to the system.  After plugging in for catalytic variables, we achieve the desired enumeration.

This theorem merits further interest because of its connection to previous work in light of the following conjecture of Duncan and Steingr\'{i}msson:

\begin{conj}[\cite{DS11}, Conjecture 3.5]
\label{DSconj}
The patterns $0021$ and $1012$ are Wilf equivalent, and $\left|\mathcal{A}_{0021}(n)\right|=\left|\mathcal{A}_{1012}(n)\right|$ is given by the binomial transform of Catalan numbers, which is sequence $A007317$ in $\cite{OEIS}$.
\end{conj}

They note that settling this conjecture would complete the Wilf classification for patterns of length 4.  Later, Mansour and Shattuck proved half of the conjecture with the following result:

\begin{theorem}[\cite{MS14}, Theorem 3.2]
\label{MSthm}
$\left|\mathcal{A}_{1012}(n)\right|=\sum_{k=0}^{n-1}\binom{n-1}{k}C_k$ for $n \geq 1$ where $C_k$ denotes the $n$th Catalan number.
\end{theorem}

The proof of this theorem is algebraic in nature making use of recurrences, nested summations, and the kernel method.

In Section \ref{fourpattern} we prove the following

\begin{theorem}
\label{BPthm2}
$\left|\mathcal{A}_{0021}(n)\right|=\sum_{k=0}^{n-1}\binom{n-1}{k}C_k$ for $n \geq 1$.
\end{theorem}

Together with Theorem \ref{MSthm}, Theorem \ref{BPthm2} answers Conjecture \ref{DSconj} in the affirmative, finishing the Wilf-classification of 4-patterns conjectured by Duncan and Steingr\'{i}msson.  The proof of Theorem \ref{BPthm2} mirrors the proof of Theorem~\ref{BPthm1} using generating trees to determine a system of multivariate generating functions from whose solution we can derive the desired enumeration.

\section{Avoiding 201 and 210} \label{pair}

\begin{customthm}{1}
$\left|\mathcal{A}_{201,210}(n)\right|=\sum_{k=0}^{n-1}\binom{n-1}{k}C_k$.
\end{customthm}

There are several components to the proof of this theorem.  They are:
\begin{enumerate}
\item Derive a generating tree for the members of $\bigcup_{n \geq 1} \mathcal{A}_{201,210}(n)$.  It turns out the nodes in our generating tree are labeled by ordered pairs $(p,q)$ where $0 \leq p < q$.  The rules for the generating tree are given in Section~\ref{gentree1}.
\item Use the generating tree from step 1 to find a recurrence for $g_{n,p,q}$, where $g_{n,p,q}$ is the number of $(p,q)$ nodes at level $n$ of the generating tree. The recurrence is given in Section \ref{rec1}.
\item Use the recurrence from step 2 to prove several relations between $g_{n,p,q}$ values for various choices of $n$, $p$, and $q$.  In particular, we set $d_{n,i} = g_{n,i-1,i}$ and $c_{n,i} = \sum_{k=0}^{i-1} g_{n,k,i}$.  The relations between the $g_{n,p,q}$ values will imply several useful relationships between the $d_{n,i}$ and $c_{n,i}$ values.  These relations are the heart of the proof of Theorem \ref{BPthm1} and are given in Section \ref{rel1}.
\item Use the relations from step 3 to derive a system of functional equations in terms of following two bivariate generating functions:
\begin{itemize} 
\item $C(x,y) = \sum_{n \geq 1} \sum_{i=1}^n c_{n,i} x^i y^n$.
\item $D(x,y) = \sum_{n \geq 1} \sum_{i=1}^n d_{n,i} x^i y^n$.
\end{itemize}
This system of equations is given in Section \ref{eq1}.
\item Although there is not a clear direct way to solve the system of functional equations in step 4, there is a unique solution.  Through computer experimentation, we conjecture the form of each of the generating functions that solve the system and verify that this set of generating functions is indeed the desired solution. The solution is given in Section~\ref{sol1}.
\item Since step 5 provides a closed form for $C(x,y) = \sum_{n \geq 1} \sum_{i=1}^n c_{n,i} x^i y^n$, we have that $C(1,y) = \sum_{n \geq 1} \sum_{i=1}^n c_{n,i}  y^n$ is the generating function for $\sum_{i=1}^n c_{n,i} = \sum_{i=1}^n \sum_{k=0}^{i-1} g_{n,k,i} = \left|\mathcal{A}_{201,210}(n)\right|$.  We then verify that $C(1,y)$ is indeed the generating function for the binomial convolution of the Catalan numbers.
\end{enumerate}

\subsection{The generating tree}\label{gentree1}

\begin{figure}
\begin{center}
\scalebox{0.8}{
\begin{tikzpicture}[level distance=1.5cm,
  level 1/.style={sibling distance=8.5cm},
  level 2/.style={sibling distance=3.5cm},
  level 3/.style={sibling distance=1cm}]
  \node {0}
    child {node {00}
      child {node {000}
      child {node {0000}}
      child {node {0001}}
      }
      child {node {001}
      child {node {0010}}
      child {node {0011}}
      child {node {0012}}
      }
    }
    child {node {01}
    child {node {010}
    child {node {0100}}
    child {node {0101}}
    child {node {0102}}
    }
      child {node {011}
      child {node {0110}}
      child {node {0111}}
      child {node {0112}}
      }
      child {node {012}
      child {node {0120}}
      child {node {0121}}
      child {node {0122}}
      child {node {0123}}
      }
    };
\end{tikzpicture}}
\end{center}
\caption{Ascent sequences of length at most 4}
\label{F01}
\end{figure}
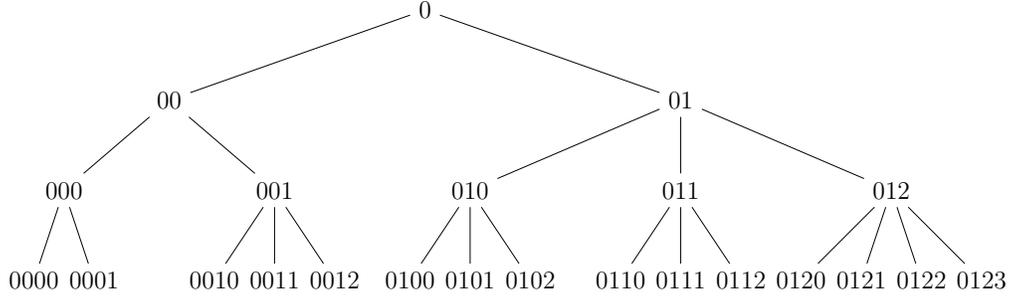

Given ascent sequences $a \in \mathcal{A}(n)$ and $a^* \in \mathcal{A}(n+1)$, we say that $a^*$ is a \emph{child} of $a$ if $a^*_1a^*_2\cdots a^*_n = a$.  In Figure \ref{F01} we see the members of $\mathcal{A}(1)$, $\mathcal{A}(2)$, $\mathcal{A}(3)$, and $\mathcal{A}(4)$ organized using the child relation.

Similarly, we may look for distinguishing features of $\{201,210\}$-avoiding ascent sequences that determine the number of $\{201,210\}$-avoiding children that a given $\{201,210\}$-avoiding ascent sequence will have.

We know that for any ascent sequence $a$, $a_i \in \{0, \dots, \mathrm{asc}(a_1\cdots a_{i-1})+1\}$.  Given $a \in \mathcal{A}_{201,210}(n)$, let $S_a$ be the set of possible integers that may be appended to $a$ to form a member of $\mathcal{A}_{201,210}(n+1)$.  For example, $S_0 = \{0,1\}$, $S_{00}=\{0,1\}$, $S_{01} = \{0,1,2\}$, and $S_{0120}=\{0,2,3\}$.

Now, consider $a \in \mathcal{A}_{201,210}(n)$ and its child $a^* \in \mathcal{A}_{201,210}(n+1)$.  Either $a^*_{n+1}>a^*_{n}$, $a^*_{n+1}=a^*_{n}$, or $a^*_{n+1}<a^*_{n}$.  We compare $S_a$ to $S_{a^*}$ in each of these three cases.

\begin{itemize}
\item If $a^*_{n+1}>a^*_{n}$, then $a^*$ has one more ascent than $a$, so $S_{a^*} = S_a \cup \{\max(S_a)+1\}$.
\item If $a^*_{n+1}=a^*_{n}$, then $S_{a^*}=S_a$.
\item If $a^*_{n+1}<a^*_{n}$, then we may not append any values between $a^*_{n}$ and $a^*_{n+1}$ (lest we form a 201 pattern), and we may not append any values smaller than $a^*_{n+1}$ (lest we form a 201 pattern).  Therefore, $S_{a^*} = S_a \setminus \left(\{i \mid a^*_{n+1} < i < a^*_{n}\} \cup \{i \mid 0 \leq i < a^*_{n+1}\}\right)$.
\end{itemize}

Notice that by definition $\left|S_a\right|$ is equal to the number of children of $a$.  Further, to determine the number of children of each child $a^*$ of $a$, we need only keep track of the last digit of $a$ and compare it to the last digit of each child $a^*$.  Therefore, the pair $(a_n, S_a)$ is sufficient to determine the pair $(a^*_{n+1}, S_{a^*})$ for each child $a^*$ of $a$.

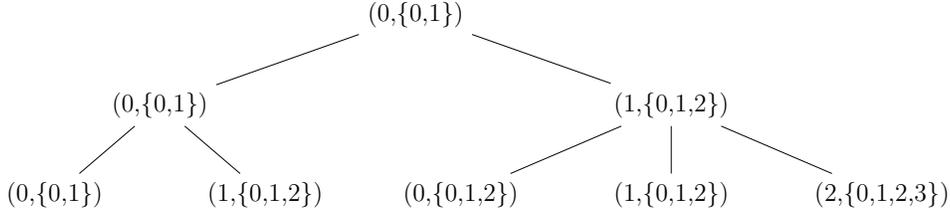
\begin{figure}
\begin{center}
\scalebox{0.8}{
\begin{tikzpicture}[level distance=1.5cm,
  level 1/.style={sibling distance=8.5cm},
  level 2/.style={sibling distance=3.5cm},
  level 3/.style={sibling distance=1cm}]
  \node {(0,\{0,1\})}
    child {node {(0,\{0,1\})}
      child {node {(0,\{0,1\})}}
      child {node {(1,\{0,1,2\})}}
    }
    child {node {(1,\{0,1,2\})}
    child {node {(0,\{0,1,2\})}}
      child {node {(1,\{0,1,2\})}}
      child {node {(2,\{0,1,2,3\})}}
    };
\end{tikzpicture}}
\end{center}
\caption{Ascent sequences $a$ of length at most 3 relabeled with the pair $(a_n, S_a)$}
\label{F02}
\end{figure}

Relabeling each ascent sequence $a$ in the first three levels of the tree in Figure \ref{F01} with the pair $(a_n, S_a)$, we obtain the tree in Figure \ref{F02}.

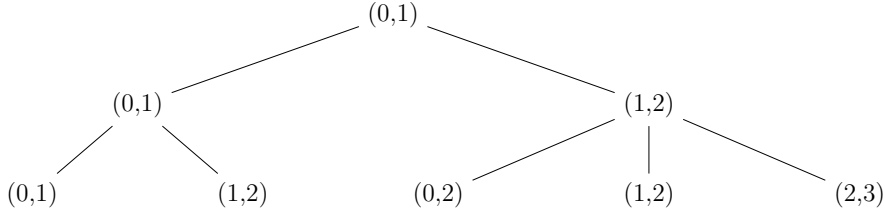
\begin{figure}
\begin{center}
\scalebox{0.8}{
\begin{tikzpicture}[level distance=1.5cm,
  level 1/.style={sibling distance=8.5cm},
  level 2/.style={sibling distance=3.5cm},
  level 3/.style={sibling distance=1cm}]
  \node {(0,1)}
    child {node {(0,1)}
      child {node {(0,1)}}
      child {node {(1,2)}}
    }
    child {node {(1,2)}
    child {node {(0,2)}}
      child {node {(1,2)}}
      child {node {(2,3)}}
    };
\end{tikzpicture}}
\end{center}
\caption{Ascent sequences $a$ of length at most 3 relabeled with the pair $(p,q)$}
\label{F03}
\end{figure}

Now, we make some normalizing conventions.  We know that $\left|S_a\right|$ gives the number of children of $a$, so it is not the particular elements of $S_a$ that matter, but rather the size of the set.  Further, we do need $a_n$ and $S_a$ to determine $S_{a^*}$ for any child $a^*$ of $a$, but again, it is not the particular digits $a_n$ and $a^*_{n+1}$ that matter, but rather how many digits of $S_a$ are smaller than, between, or larger than these digits.  Therefore, let the \emph{reduction} of the pair $(a_n, S_a)$ be the pair $\mathrm{red}((a_n,S_a))$ obtained by replacing the $i$th smallest digits with $i-1$.  For example, $\mathrm{red}((5,\{1,4,5,7\})) = (2,\{0,1,2,3\})$.  Notice that $a_n \in S_a$ for all ascent sequences $a$ since repeating the last digit of a given $\{201,210\}$-avoiding ascent sequence still produces an ascent sequence that avoids 201 and 210.  This means that if $(a_n^{\prime},S_a^{\prime})=\mathrm{red}((a_n,S_a))$, then $S_a^{\prime}$ is a set of consecutive integers with minimum 0 and $a_n^{\prime} \in S_a^{\prime}$.  Therefore, we may more concisely represent $\mathrm{red}((a_n,S_a))$ by the ordered pair $(p,q)$ where $p$ is the first element of $\mathrm{red}((a_n,S_a))$ and $q$ is the maximum of the second element of $\mathrm{red}((a_n,S_a))$.    To this end, the tree from Figure \ref{F02} may be relabeled as seen in Figure \ref{F03} using reduction.  As a more interesting example later in the generating tree, $a = 0101341 \in \mathcal{A}_{201,210}(7)$ has $(a_n, S_a) = (1,\{1,4,5\})$, but after reduction, we obtain $\mathrm{red}((a_n,S_a)) = (0,\{0,1,2\})$, which we relabel again as $(0,2)$.

Using this labeling and the rules for computing $S_{a^*}$ from $S_a$ given above, we build a generating tree for $\mathcal{A}_{201,210}$.  The root $(0,1)$ comes from the definitions and labeling conventions above.  Given $(p,q)$, which is shorthand for $a_n=p$, $S_a=\{0,\dots, q\}$, we see that if we append a digit $i$ from $\{p+1,\dots, q\}$ to $a$, then $S_{a^*} = \{0,\dots, q+1\}$, so we obtain an ascent sequence with label $(i,q+1)$.  If we append another copy of $p$, we obtain another ascent sequence with label $(p,q)$. If we append a digit $i$ from $\{0, \dots, p-1\}$, we obtain an ascent sequence $a^*$ with $S_{a^*} = \{i,p,p+1,\dots, q\}$; since $\mathrm{red}((i,\{i,p,p+1,\dots, q\})) = (0,\{0,\dots, 1+q-p\})$ this node has label $(0,1+q-p)$.

In particular, the root and rules for our generating tree are
\begin{itemize}
\item \textbf{root:} $(0,1)$
\item \textbf{rule:} $(p,q) \to (0,1+q-p)^p, (p,q), (p+1,q+1), (p+2,q+1), (p+3,q+1), \dots, (q,q+1)$
\end{itemize}
where $(0,1+q-p)^p$ indicates that we have $p$ copies of the node $(0,1+q-p)^p$.

Notice that by convention, $0 \leq p<q$ for all nodes in this generating tree.  Further, $q$ increases by at most 1 with each generation of the tree, so the largest value of $q$ in a node at level $n$ is $q=n$.  We wish to find a formula for the number of nodes at level $n$ of this generating tree.

\subsection{The recurrence} \label{rec1}

Now, let $g_{n,p,q}$ be the number of nodes of type $(p,q)$ at level $n$ in the generating tree of Section \ref{gentree1}.  Certainly, $g_{n,p,q}=0$ if $n \leq 0$.

We saw above that $g_{n,p,q} = 0$ if $p \geq q$, if $n=1$ and $p \neq 0$, if $n=1$ and $q \neq 1$, or if $q \geq n+1$.  We also know that $g_{1,0,1}=1$, which counts the root node of the generating tree.

We still must determine $g_{n,p,q}$ for other values.  We consider 2 cases: when $p=0$ and when $p>0$.

In the case where $p=0$, we obtain a $(0,q)$ node in level $n$ for every $(0,q)$ node at level $n-1$, and we obtain $i$ $(0,q)$ nodes at level $n$ for every $(i,q-1+i)$ node at level $n$.  Therefore $g_{n,0,q} = g_{n-1,0,q} + \sum_{i=1}^{n-q} i\cdot g_{n-1,i,q-1+i}$.

Finally, in the case where $p>0$, we obtain a $(p,q)$ node at level $n$ for every $(p,q)$ node at level $n-1$, and we obtain a $(p,q)$ node at level $n$ for every $(i,q-1)$ node at level $n-1$ where $0 \leq i \leq p-1$.  Therefore $g_{n,p,q} = g_{n-1,p,q} + \sum_{i=0}^{p-1} g_{n-1,i,q-1}$.

To summarize, $\left|\mathcal{A}_{201,210}(n)\right| = \sum_{q=1}^n \sum_{p=0}^{q-1} g_{n,p,q}$ where

\begin{equation}
g_{n,p,q} = \begin{cases}
0&n \leq 0\\
0& p \geq q\\
0& n=1 \text{ and } p \neq 0\\
0& n=1 \text{ and } q \neq 1\\
1 & n=1 \text{ and } p=0 \text{ and } q=1\\
g_{n-1,0,q} + \sum\limits_{i=1}^{n-q}i \cdot g_{n-1,i,q-1+i}&p=0\\
g_{n-1,p,q} + \sum\limits_{i=0}^{p-1}g_{n-1,i,q-1}& \text{ otherwise.}
\end{cases}
\label{Recur1}
\end{equation}

\subsection{Relations} \label{rel1}

We have that $\left|\mathcal{A}_{201,210}(n)\right| = \sum_{q=1}^n \sum_{p=0}^{q-1} g_{n,p,q}$.

We wish to understand the structure of $g_{n,p,q}$ more fully.  To this end, let
\begin{itemize}
\item $c_{n,i} = \sum_{k=0}^{i-1} g_{n,k,i}$ (i.e. $c_{n,i}$ is the number of nodes of type $(*,i)$ at level $n$ in the generating tree.)
\item $d_{n,i} = g_{n,i-1,i}$.
\end{itemize}
Notice that $\left|\mathcal{A}_{201,210}(n)\right|=\sum_{i=1}^n c_{n,i}$ by definition.

For fixed $n$, let $A_n$ be the $n \times n$ array containing $g_{n,p,q}$ in row $p+1$, column $q$.  In this arrangement, $d_{n,i}$ is the $i$th entry on the main diagonal of $A_n$ and $c_{n,i}$ is the sum of the entries in the $i$th column.  Arrays $A_j$ ($1 \leq j \leq 7$) and corresponding values of $c_{n,i}$ and $d_{n,i}$ are given in Table \ref{arrays1}.  Organizing the numbers $g_{n,p,q}$ in this two dimensional format reveals some interesting patterns.

\begin{table}
\begin{center}
\begin{tabular}{|l|l|l|}
\hline
$A_1 = \left[\begin{array}{c}1\end{array} \right]$&$c_{1,1}=1$&$d_{1,1}=1$\\
\hline
\multirow{2}{*}{$A_2 = \left[\begin{array}{cc}1&0\\0&1\end{array} \right]$}&$c_{2,1}=1$&$d_{2,1}=1$\\
&$c_{2,2}=1$&$d_{2,2}=1$\\
\hline
\multirow{3}{*}{$A_3 = \left[\begin{array}{ccc}1&1&0\\0&2&0\\0&0&1\end{array} \right]$}&$c_{3,1}=1$&$d_{3,1}=1$\\
&$c_{3,2}=3$&$d_{3,2}=2$\\
&$c_{3,3}=1$&$d_{3,3}=1$\\
\hline
\multirow{4}{*}{$A_4 = \left[\begin{array}{cccc}1&5&0&0\\0&3&1&0\\0&0&4&0\\0&0&0&1\end{array} \right]$}&$c_{4,1}=1$&$d_{4,1}=1$\\
&$c_{4,2}=8$&$d_{4,2}=3$\\
&$c_{4,3}=5$&$d_{4,3}=4$\\
&$c_{4,4}=1$&$d_{4,4}=1$\\
\hline
\multirow{5}{*}{$A_5 = \left[\begin{array}{ccccc}1&19&1&0&0\\0&4&6&0&0\\0&0&12&1&0\\0&0&0&6&0\\0&0&0&0&1\end{array} \right]$}&$c_{5,1}=1$&$d_{5,1}=1$\\
&$c_{5,2}=23$&$d_{5,2}=4$\\
&$c_{5,3}=19$&$d_{5,3}=12$\\
&$c_{5,4}=7$&$d_{5,4}=6$\\
&$c_{5,5}=1$&$d_{5,5}=1$\\
\hline
\multirow{6}{*}{$A_6 = \left[\begin{array}{cccccc}1&69&9&0&0&0\\0&5&25&1&0&0\\0&0&35&8&0&0\\0&0&0&25&1&0\\0&0&0&0&8&0\\0&0&0&0&0&1\end{array} \right]$}&$c_{6,1}=1$&$d_{6,1}=1$\\
&$c_{6,2}=74$&$d_{6,2}=5$\\
&$c_{6,3}=69$&$d_{6,3}=35$\\
&$c_{6,4}=34$&$d_{6,4}=25$\\
&$c_{6,5}=9$&$d_{6,5}=8$\\
&$c_{6,6}=1$&$d_{6,6}=1$\\
\hline
\multirow{7}{*}{$A_7 = \left[\begin{array}{ccccccc}1&256&53&1&0&0&0\\0&6&94&10&0&0&0\\0&0&109&42&1&0&0\\0&0&0&94&10&0&0\\0&0&0&0&42&1&0\\0&0&0&0&0&10&0\\0&0&0&0&0&0&1\end{array} \right]$}&$c_{7,1}=1$&$d_{7,1}=1$\\
&$c_{7,2}=262$&$d_{7,2}=6$\\
&$c_{7,3}=256$&$d_{7,3}=109$\\
&$c_{7,4}=147$&$d_{7,4}=94$\\
&$c_{7,5}=53$&$d_{7,5}=42$\\
&$c_{7,6}=11$&$d_{7,6}=10$\\
&$c_{7,7}=1$&$d_{7,7}=1$\\
\hline
\end{tabular}
\end{center}
\caption{Arrays containing the values of $g_{n,p,q}$ for $n \leq 7$}
\label{arrays1}
\end{table}

\begin{lemma}\label{BigLemma}
The following relations hold for $g_{n,p,q}$, $d_{n,i}$, and $c_{n,i}$:
   \begin{enumerate}[label={\thelemma.\alph*}]
       \item $d_{n,i}=d_{n-1,i}+c_{n-1,i-1}$ for $2 \leq i \leq n$. \label{L7}
       \item $g_{n,i,n}=0$ for $n \geq 1$ and $i \neq n-1$. (The rightmost column consists of $0$s except for the bottom entry.) \label{L5} 
			 \item $g_{n,0,1}=d_{n,1}=c_{n,1}=d_{n,n}=c_{n,n}=1$ for $n \geq 1$. (The top left entry and the bottom right entry of $A_n$ are both $1$, but the rest of the first column and the rest of the last column are all $0$s.) \label{L6} 
			 \item $g_{n,p,q} = g_{n,p+2,q+1}$ for $3 \leq p+2 < q+1 \leq n$. (This implies $p \geq 1$ and $q \geq 3$). (We have a recursive way to compute non-diagonal entries other than in the first row of $A_n$.) \label{L8} 
			 \item $g_{n,0,q} = g_{n,0,q+1} + g_{n,1,q+1} + g_{n,2,q+1}$ for $2 \leq q \leq n-1$. (We have a recursive way to compute the entries of the first row of $A_n$ from other entries.) \label{L9}
			 \item $\sum_{i=0}^k g_{n,i,q} = \sum_{i=0}^{k+2} g_{n,i,q+1}$ for $k\leq q-2$. (This gives a condition on partial sums of adjacent columns in $A_n$.) \label{L10}
			 \item $c_{n,i}=c_{n,i-1}-d_{n,i-1}$ for $3 \leq i \leq n$. \label{L11}
   \end{enumerate}
\end{lemma}

In particular, part \ref{L7} follows directly from Equation \ref{Recur1}.  Parts \ref{L5} and \ref{L6} are used in the proof of parts \ref{L8} and \ref{L9}.  Parts \ref{L8} and \ref{L9} are used to prove part \ref{L10}, which is used for part \ref{L11}.  Notice that given array $A_{n-1}$, parts \ref{L7} and \ref{L6} completely determine the diagonal entries of $A_{n}$, and then parts \ref{L6} and \ref{L11} completely determine the column sums from the diagonals (working from right to left).  Thus, only parts \ref{L7}, \ref{L6}, and \ref{L11} are used to derive functional equations in Section \ref{eq1}.  While it is easy to verify that Lemma \ref{BigLemma} holds for the arrays given in Table \ref{arrays1}, the proof for the general case is long and technical.  The interested reader can find details of the proof in Appendix \ref{Appendix}.

\subsection{Functional equations} \label{eq1}

From the previous section we know:

\begin{itemize}
\item $c_{n,i} = c_{n,i-1} - d_{n,i-1}$ for $3 \leq i \leq n$,
\item $d_{n,i} = d_{n-1,i} + c_{n-1,i-1}$ for $2 \leq i \leq n$.
\end{itemize}

Now, define the following three generating functions:

\begin{itemize}
\item $C(x,y) = \sum_{n \geq 1} \sum_{i=1}^n c_{n,i} x^iy^n$,
\item $D(x,y) = \sum_{n \geq 1} \sum_{i=1}^n d_{n,i} x^iy^n$,
\item $C_2(y) = \sum_{n \geq 2} c_{n,2} y^n$.
\end{itemize}

$c_{n,i} = c_{n,i-1} - d_{n,i-1}$ for $3 \leq i \leq n$ implies that 

\begin{equation}
(1-x)C(x,y) + xD(x,y) = \frac{xy}{1-y} + x^2C_2(y).
\end{equation}

$d_{n,i} = d_{n-1,i} + c_{n-1,i-1}$ for $2 \leq i \leq n$ implies that

\begin{equation}
(1-y)D(x,y) - xyC(x,y) = xy.
\end{equation}

We also know that $c_{n,n}=1$ for $n \geq 1$, which implies that 

\begin{equation}
\begin{split}
C\left(\frac{1}{y},yz\right)\Biggr|_{y=0} & = \left(\sum_{n \geq 1} \sum_{i=1}^n c_{n,i} \left(\frac{1}{y}\right)^i(yz)^n\right)\Biggr|_{y=0} \\
 & = \left(\sum_{n \geq 1} \sum_{i=1}^n c_{n,i} y^{n-i}z^n\right)\Biggr|_{y=0}\\
 & = \sum_{n \geq 1} \sum_{i=1}^n c_{n,i} 0^{n-i}z^n \\
 & = \sum_{n \geq 1} c_{n,n}z^n \\
 & = \frac{z}{1-z}.
\end{split}
\end{equation}

\subsection{Generating functions} \label{sol1}

Now, we wish to solve the following system of equations.

\begin{itemize}
\item $(1-x)C(x,y) + xD(x,y) = \frac{xy}{1-y} + x^2C_2(y)$,
\item $(1-y)D(x,y) - xyC(x,y) = xy$,
\item $C\left(\frac{1}{y},yz\right)\Biggr|_{y=0} = \dfrac{z}{1-z}$.
\end{itemize}

Note that it was necessary to introduce $C_2(y)$ separately since the recurrence for $c_{n,i}$ only applies for $i \geq 3$.  The first two equations are linear in $C(x,y)$ and $D(x,y)$, while the third equation puts a condition on the coefficient of $x^ny^n$ in $C(x,y)$. 

It turns out that there are infinitely many solutions to the first two equations, but the fact that $c_{n,n}=1$ for $n \geq 1$ determines the unique solution.  To be sure:

Multiply the first equation by $(1-y)$ and the second equation by $x$ to obtain
\begin{itemize}
\item $(1-x)(1-y)C(x,y) + x(1-y)D(x,y) = xy + x^2(1-y)C_2(y)$,
\item $x(1-y)D(x,y) - x^2yC(x,y) = x^2y$.
\end{itemize}
After subtracting the second equation from the first, we have:
\begin{equation}
\left((1-x)(1-y)+x^2y\right)C(x,y) = xy-x^2y+x^2(1-y)C_2(y),
\end{equation}
or, equivalently,
\begin{equation}
C(x,y) = \frac{xy-x^2y+x^2(1-y)C_2(y)}{(1-x)(1-y)+x^2y}.
\end{equation}

If we replace $C_2(y)$ with the formal power series $\sum_{n \geq 2} c_{n,2} y^n$, then, after expanding, we see that the coefficient of $xy$ in $C(x,y)$ is 1, and the coefficient of $x^iy^i$ in $C(x,y)$ is a linear expression in terms of $c_{2,2}, \dots, c_{i,2}$.  Lemma \ref{BigLemma} shows that $c_{n,n}=1$ for all $n \geq 1$, so the coefficient of $x^iy^i$ in $C(x,y)$ is 1.  This implies a unique solution for the values of $c_{2,2}, \dots, c_{i,2}$.

This additional fact implies that there is a unique set of three functions that satisfy the system of equations where $c_{n,n}=1$ for all $n$, but there is not a straightforward way to solve for the functions directly.  However, using the \texttt{gfun} package in Maple, we can predict the form of $C_2(y)$.  If we know $C_2(y)$, we can plug it into the first equation and then use the first two equations to find conjectured forms for $C(x,y)$ and $D(x,y)$.

It turns out that $C(x,y)$, $D(x,y)$, and $C_2(y)$ are as follows:

\begin{equation}
C(x,y) = \frac{\left(x\sqrt{5y^2-6y+1}-xy+x+2y-2\right)xy}{2(x^2y+xy-x-y+1)(y-1)},
\label{finalC}
\end{equation}

\begin{equation}
\begin{aligned}
&D(x,y) =\\
& \frac{-\left(x^2y\sqrt{5y^2-6y+1}+x^2y^2-x^2y+4xy^2-6xy-2y^2+2x+4y-2\right)xy}{2(x^2y^2-x^2y+xy^2-2xy-y^2+x+2y-1)(y-1)},
\end{aligned}
\end{equation}

\begin{equation}
C_2(y) = \frac{-y\left(-1+y+\sqrt{5y^2-6y+1}\right)}{2(y-1)^2}.
\end{equation}

It is straightforward to plug these three equations into the original system of two equations and verify that they are \emph{a} solution.  It can also be checked that for the expression $C(x,y)$ given in Equation \ref{finalC} we have $C(\frac{1}{y},yz)\bigr|_{y=0}=\frac{z}{1-z}$, which implies $c_{n,n}=1$ for all $n \geq 1$.

\subsection{The punchline}

We have determined a closed form for $C(x,y) = \sum_{n \geq 1} \sum_{i=1}^n c_{n,i} x^i y^n$, so we have that $C(1,y) = \sum_{n \geq 1} \sum_{i=1}^n c_{n,i}  y^n$ is the generating function for $\sum_{i=1}^n c_{n,i} = \sum_{i=1}^n \sum_{k=0}^{i-1} g_{n,k,i} = \left|\mathcal{A}_{201,210}(n)\right|$.
\begin{equation}
C(1,y) = \frac{-1+y+\sqrt{5y^2-6y+1}}{2(y-1)}
\end{equation}
This is, as per OEIS entry A007317, the generating function for the binomial convolution of the Catalan numbers.

We have now seen that $\{201,210\}$-avoiding ascent sequences are enumerated by the binomial convolution of the Catalan numbers.  To be sure, our proof show inherent structure in the set $\mathcal{A}_{201,210}(n)$, but it requires an experimental prediction that is later validated.  It remains open to find a statistic $\mathrm{st}: \mathcal{A}_n \to \mathbb{N}$ such that $\left|\left\{a\in \mathcal{A}_{201,210}(n) \mid \mathrm{st}(a)=k\right\}\right| = \binom{n-1}{k}C_k$.

\section{Avoiding 0021} \label{fourpattern}

Our experiment-based methodology in Section \ref{pair} is sufficiently general that we next adapt it to enumerate 0021-avoiding ascent sequences.  The enumeration of such sequences was also conjectured to be the binomial convolution of the Catalan numbers in \cite{DS11}, but a proof has remained open until now.  Combining this work with previous results of Duncan and Steingr\'{i}msson \cite{DS11}, and Mansour and Shattuck \cite{MS14}, we complete the Wilf-classification of 4-patterns in the context of ascent sequences.

\begin{customthm}{4}
$\left|\mathcal{A}_{0021}(n)\right|=\sum_{k=0}^{n-1}\binom{n-1}{k}C_k$.
\end{customthm}

Again there are several components to the proof of this theorem.  The enumeration argument mirrors our approach to Theorem \ref{BPthm1}.  While the generating tree is more complicated, less work is required to convert the generating tree rules into a system of functional equations.  The outline of the proof is given below.
\begin{enumerate}
\item Derive a generating tree for the members of $\bigcup_{n \geq 1} \mathcal{A}_{0021}(n)$.  It turns out the nodes in our generating tree are labeled by ordered triples $(p,q,r)$ where $p,q,r \geq 0$ and $p \in \{q-2,q-1,q\}$.  This generating tree is shown in Section \ref{gentree2}.
\item Use the generating tree from step 1 to find a recurrence for $g_{n,p,q,r}$, where $g_{n,p,q,r}$ is the number of $(p,q,r)$ nodes at level $n$ of the generating tree.  Since $p \in \{q-2,q-1,q\}$, we will consider $g_{n,q-2,q,r}$, $g_{n,q-1,q,r}$ and $g_{n,q,q,r}$ values separately. This analysis is given in Section \ref{rec2}.
\item Use the recurrence from step 3 to derive a system of two functional equations in terms of two following trivariate generating functions:
\begin{itemize} 
\item $C(x,y,z) = \sum_{n \geq 1} \sum_{q \geq 1} \sum_{r \geq 2} g_{n,q,q,r} x^q y^r z^n$,
\item $D(x,y,z) = \sum_{n \geq 1} \sum_{q \geq 1} \sum_{r \geq 1} g_{n,q-1,q,r} x^q y^r z^n$.
\end{itemize}
\item Although there is not a clear direct way to solve the system of functional equations in step 3, since there are two equations and two unknown functions, there must be a unique solution.  Through computer experimentation, we conjecture the form of each of the generating functions that solve the system and verify that this set of generating functions is indeed the desired solution.  The functional equations from step 3 and their solution are given in Section \ref{eq2}.
\item Since step 4 provides closed forms for $$C(x,y,z)= \sum_{n \geq 1} \sum_{q \geq 1} \sum_{r \geq 2} g_{n,q,q,r} x^q y^r z^n$$ and $$D(x,y,z) = \sum_{n \geq 1} \sum_{q \geq 1} \sum_{r \geq 1} g_{n,q-1,q,r} x^q y^r z^n,$$ and there is exactly one node of type $(q-2,q,r)$ on each level, we have that $C(1,1,z)+D(1,1,z)+\frac{z}{1-z}$ is the generating function for $\left|\mathcal{A}_{0021}(n)\right|$.  We then verify that $C(1,1,z)+D(1,1,z)+\frac{z}{1-z}$  is indeed the generating function for the binomial convolution of the Catalan numbers.
\end{enumerate}

\subsection{The generating tree} \label{gentree2}

As in Section \ref{pair}, we first organize the members of $\bigcup_{n \geq 1} \mathcal{A}(n)$ using the child relation.  The resulting tree was given in Figure \ref{F01}.  Next, we apply the child relation to the members of $\bigcup_{n \geq 1} \mathcal{A}_{0021}(n)$ and look for distinguishing features of $\{0021\}$-avoiding ascent sequences that determine the number of children a given ascent sequence will have.

We still know that for any ascent sequence $a$, $a_i \in \{0, \dots, \mathrm{asc}(a_1\cdots a_{i-1})+1\}$.  Given $a \in \mathcal{A}_{0021}(n)$, let $S_a$ be the set of possible integers that may be appended to $a$ to form a member of $\mathcal{A}_{0021}(n+1)$.  For example, $S_0 = \{0,1\}$, $S_{00}=\{0,1\}$, $S_{01} = \{0,1,2\}$, and $S_{01013}=\{0,3,4\}$.

Notice that members of $S_a$ are affected by repeated digits appearing in $a$.  In particular, if $a_i=a_j=x$ where $i <j$ then all digits after position $j$ and larger than $x$ must appear in decreasing order.  Also, if $a_i=a_j=x$ and $a_k=y>x$ where $i<j<k$, then no digits from $\{x+1,\dots,y-1\}$ may be appended in the rest of the ascent sequence.  Therefore, we break $S_a$ into $S_a^u$ and $S_a^i$.  $S_a^i =\{x \in S_a \mid x > \text{(smallest repeated digit in $a$)}\}$, and $S_a^u = S_a \setminus S_a^i$.  Here, the superscripts $u$ and $i$ stand for ``unrestricted'' vs. ``increasing'' respectively since all digits larger than and after a repeated digit must appear in increasing order.  For example, if $a=0121235$, then $S_a = \{0,1,5,6\}$.  Since the smallest repeated digit in $a$ is 1, we have $S_a^u = \{0,1\}$ and $S_a^i=\{5,6\}$. Alternatively, if $a=01234$ then $S_a=\{0,1,2,3,4,5\}$, and since there is no smallest repeated digit, $S_a^u=\{0,1,2,3,4,5\}$ while $S_a^i = \emptyset$.

Now, consider $a \in \mathcal{A}_{0021}(n)$ and its child $a^* \in \mathcal{A}_{0021}(n+1)$.  We compare $S_a^u$ and $S_a^i$ to $S_{a^*}^u$ and $S_{a^*}^i$ in each of several cases.

\begin{itemize}
\item If $a_{n+1}^* \leq a_n$, then we have not created a new ascent so $S_a^u \cup S_a^i = S_{a^*}^u \cup S_{a^*}^i$, but we may have created a new smallest repeated digit.  If $a_{n+1}^* \in \{\max(S_a^u)\} \cup S_a^i$, then we have not created a new smallest repeated digit.  If $a_{n+1}^* \in S_a^u \setminus \{\max(S_a^u)\}$, then we have created a new smallest repeated digit.  We modify $S_a^u$ and $S_a^i$ as follows:
\begin{itemize}
\item If $a_{n+1}^* \leq a_n$ and $a_{n+1}^* \in \{\max(S_a^u)\} \cup S_a^i$, then $S_{a^*}^u = S_a^u$ and $S_{a^*}^i = S_a^i$.

For example, if $a^*=01212351$ or $a^*=01212355$, then $S_{a^*}^u=\{0,1\}$ and $S_{a^*}^i=\{5,6\}$.
\item If $a_{n+1}^* \leq a_n$ and $a_{n+1}^* \in S_a^u \setminus \{\max(S_a^u)\}$, then all digits larger than $a_{n+1}^*$ must be moved to $S_{a^*}^i$, so $S_{a^*}^u = \{j \in S_a^u \mid j \leq a_{n+1}^*\}$ and $S_{a^*}^i = S_a^i \cup \{j \in S_a^u \mid j > a_{n+1}^*\}$.

For example, if $a^*=01212350$, then $S_{a^*}^u=\{0\}$ and $S_{a^*}^i=\{1,5,6\}$. 
\end{itemize}
\item If $a_{n+1}^* > a_n$, then we have created a new ascent, so $\max(S_{a^*}^u \cup S_{a^*}^i) = 1+\max(S_a^u \cup S_a^i)$.
\begin{itemize}
\item If $a_{n+1}^* \in S_a^u$, then $a$ must be the strictly increasing sequence of length $n$, so $S_{a^*}^i=S_a^i = \emptyset$, while $S_{a^*}^u = S_a^u \cup \{a_{n+1}^*+1\}$.

For example, if $a^*=012345$, then $S_{a^*}=S_{a^*}^u = \{0,1,2,3,4,5,6\}$ and $S_{a^*}^i=\emptyset$.
\item If $a_{n+1}^* \in S_a^i$ then $S_{a^*}^u = S_{a}^u$ and $S_{a^*}^i = S_a^i \cup \{\max(S_a^i)+1\} \setminus \{j \in S_a^i \mid j < a_{n+1}^*\}$.

For example, if $a^*=01212356$, then $S_{a^*}^u=\{0,1\}$ and $S_{a^*}^i=\{6,7\}$.
\end{itemize}
\end{itemize}

Notice that by definition $\left|S_a^u \cup S_a^i\right|$ is equal to the number of children of $a$.  Further, to determine the number of children of each child $a^*$ of $a$, we need only keep track of the last digit of $a$ and compare it to the last digit of each child $a^*$.  Therefore, the triple $(a_n, S_a^u, S_a^i)$ is sufficient to determine the triple $(a^*_{n+1}, S_{a^*}^u, S_{a^*}^i)$ for each child $a^*$ of $a$.

\begin{figure}
\begin{center}
\scalebox{0.8}{
\begin{tikzpicture}[level distance=1.5cm,
  level 1/.style={sibling distance=8.5cm},
  level 2/.style={sibling distance=3.5cm},
  level 3/.style={sibling distance=1cm}]
  \node {(0,\{0,1\},\{\})}
    child {node {(0,\{0\},\{1\})}
      child {node {(0,\{0\},\{1\})}}
      child {node {(1,\{0\},\{1,2\})}}
    }
    child {node {(1,\{0,1,2\},\{\})}
    child {node {(0,\{0\},\{1,2\})}}
      child {node {(1,\{0,1\},\{2\})}}
      child {node {(2,\{0,1,2,3\},\{\})}}
    };
\end{tikzpicture}}
\end{center}
\caption{Ascent sequences $a$ of length at most 3 relabeled with the triple $(a_n, S_a^u, S_a^i)$}
\label{F02b}
\end{figure}
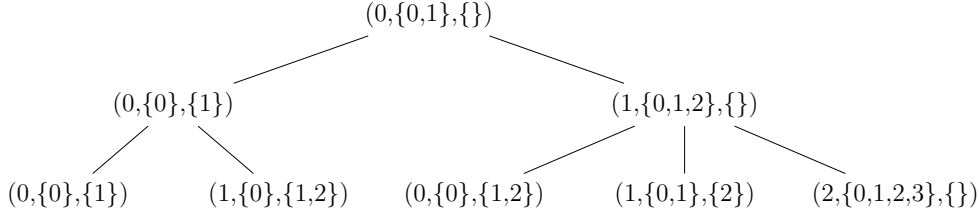

Relabeling each ascent sequence $a$ in the first three levels of the tree in Figure \ref{F01} with the triple $(a_n, S_a^u, S_a^i)$, we obtain the tree in Figure \ref{F02b}.

Now, we make some normalizing conventions.  We know that $\left|S_a\right| = \left| S_a^u \cup S_a^i\right|$ gives the number of children of $a$, so it is not the particular elements of $S_a$ that matter, but rather the size of the set.  Further, we do need $a_n$ and $S_a$ to determine $S_{a^*}$ for any child $a^*$ of $a$, but again, it is not the particular digits $a_n$ and $a^*_{n+1}$ that matter, but rather how many digits of $S_a$ are smaller than, between, or larger than these digits.  Therefore, let the \emph{reduction} of the triple $(a_n, S_a^u, S_a^i)$ be the triple $\mathrm{red}((a_n,S_a^u,S_a^i))$ obtained by replacing the $i$th smallest digits with $i-1$.  For example, $\mathrm{red}((2,\{0\},\{2,3\})) = (1,\{0\},\{1,2\})$.  

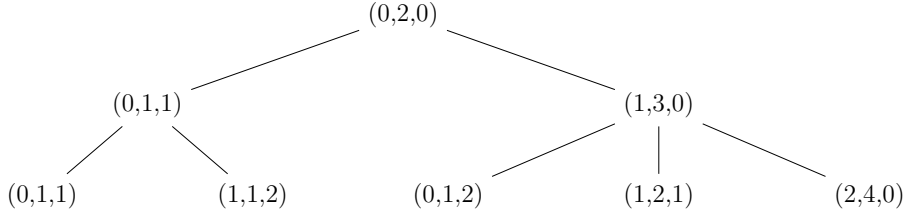
\begin{figure}
\begin{center}
\scalebox{0.8}{
\begin{tikzpicture}[level distance=1.5cm,
  level 1/.style={sibling distance=8.5cm},
  level 2/.style={sibling distance=3.5cm},
  level 3/.style={sibling distance=1cm}]
  \node {(0,2,0)}
    child {node {(0,1,1)}
      child {node {(0,1,1)}}
      child {node {(1,1,2)}}
    }
    child {node {(1,3,0)}
    child {node {(0,1,2)}}
      child {node {(1,2,1)}}
      child {node {(2,4,0)}}
    };
\end{tikzpicture}}
\end{center}
\caption{Ascent sequences $a$ of length at most 3 relabeled with the triple $(p,q,r)$}
\label{F03b}
\end{figure}

Notice that $a_n \in S_a$ for all ascent sequences $a$ since repeating the last digit of a given $\{0021\}$-avoiding ascent sequence still produces an ascent sequence that avoids 0021.  This means that if $(a_n^\prime, S_a^{u^\prime}, S_a^{i^\prime})= \mathrm{red}((a_n,S_a^u,S_a^i))$ then $S_a^{u^\prime} \cup S_a^{i^\prime}$ is a set of consecutive integers with minimum 0 and $a_n^\prime \in S_a^{u^\prime} \cup S_a^{i^\prime}$.  Therefore, we may more concisely represent $\mathrm{red}((a_n,S_a^u, S_a^i))$ by the ordered triple $(p,q,r)$ where $p$ is the first element of $\mathrm{red}((a_n,S_a^u,S_a^i))$ and $q=\left|S_a^u\right|$ and $r=\left|S_a^i\right|$.    To this end, the tree from Figure \ref{F02b} may be relabeled as seen in Figure \ref{F03b}.  As a more interesting example, $0102 \in \mathcal{A}_{0021}(4)$ would initially be relabeled as $(2,\{0\},\{2,3\})$, but after reduction, we obtain $(1,\{0\},\{1,2\})$, which we relabel again as $(1,1,2)$.

Notice that for a node of type $(p,q,r)$, we have $p \in \{q-2,q-1,q\}$.  That is, $a_n$ is always one of the largest two values in $S_a^u$ or is it the smallest value in $S_a^i$ after reduction.

Using this labeling and the rules for computing $S_{a^*}^u$ and $S_{a^*}^i$ from $S_a^u$ and $S_a^i$ given above, we may build a generating tree isomorphic to the generating tree for $\mathcal{A}_{0021}$.  The root $(0,2,0)$ comes from the definitions and labeling conventions above.

A node of type $(q-2,q,0)$ is a strictly increasing ascent sequence of length $q-1$ and has $q$ children, one from appending each of the digits $\{0,1,\dots, q-1\}$.  If we append $q-1$, we create another increasing sequence with label $(q-1, q+1, 0)$.  If we append a smaller digit $i$, then we have repeated that digit, so $S_{a^*}^u = \{0, \dots, i\}$, and $S_{a^*}^i = \{i+1, \dots, q-1\}$.

A node of type $(q-1,q,r)$ has $q+r$ children.  If we repeat the last digit, $S_a^u$ and $S_a^i$ remain unchanged.  If we append one of the $q-2$ digits $i \in S_a^u$ that is smaller than $q-1$, we have repeated that digit so $S_{a^*}^u = \{0, \dots, i\}$, and $S_{a^*}^i = \{i+1, \dots, q+r-1\}$.  If we append one of the $r$ digits $j$ of $S_u^i$, then $j$ becomes the new smallest digit of $S_a^i$ while the size of $S_a^i$ increases by 1, forming a node of type $(q,q,j)$.

A node of type $(q,q,r)$ has $q+r$ children.  If we repeat the last digit, $S_a^u$ and $S_a^i$ remain unchanged.  If we append one of the $q$ digits $i \in S_a^u$ that is smaller than $q$, $S_{a^*}^u = \{0,\dots, i\}$ and $S_a^i$ decreases accordingly.  If we append one of the $r-1$ digits $j$ of $S_a^i$ that is larger than $q$, we obtain a node of type $(q,q,j)$.

In particular we have
\begin{itemize}
\item \textbf{root:} $(0,2,0)$
\item \textbf{rules:}

$(q-2,q,0) \to (q-1,q+1,0), (i,i+1,q-1-i)_{i=0}^{q-2}$

$(q-1,q,r) \to (q-1,q,r), (i,i+1,q+r-1-i)_{i=1}^{q-2}, (q,q,i)_{i=2}^{r+1}$

$(q,q,r) \to (q,q,r), (i,i+1,q+r-1-i)_{i=0}^{q-1}, (q,q,i)_{i=2}^{r}$

\end{itemize}
where $(a(i),b(i),c(i))_{i=d}^e$ corresponds to a list of $e-d+1$ nodes where we plug in each value $i \in \{d,d+1,\dots, e\}$ into the expressions $a(i)$, $b(i)$, and $c(i)$.

Notice that by convention all nodes are of the form $(p,q,r)$ where $p,q,r \in \mathbb{N}$, $p \in \{q-2,q-1,q\}$ and 
\begin{itemize}
\item If $p=q-2$, then $r=0$.
\item If $p=q-1$, then $q \geq 1$, $r \geq 1$, and $q+r \leq n$.
\item If $p=q$, then $q \geq 1$, $r \geq 2$, $q+r \leq n$.
\end{itemize}

We wish to find a formula for the number of nodes at level $n$ of this generating tree.

\subsection{The recurrence} \label{rec2}

Now, let $g_{n,p,q,r}$ be the number of nodes of type $(p,q,r)$ at level $n$ in the generating tree above.  Certainly, $g_{n,p,q,r}=0$ if $n \leq 0$.

As above, we know that $p \in \{q-2,q-1,q\}$, so let $g^0_{n,q,r} = g_{n,q,q,r}$, $g^1_{n,q,r} = g_{n,q-1,q,r}$, and $g^2_{n,q,r} = g_{n,q-2,q,r}$.  We have the following:

\begin{equation}
g^2_{n,q,r} = 0 \text{ if } r \neq 0,
\end{equation}

\begin{equation}
g^2_{n,q,0} = g^2_{n-1,q-1,0},
\end{equation}

\begin{equation}
\begin{split}
&g^0_{n,q,r} =\\
& \quad \begin{cases}
0&q=0\\
0& n \leq 2\\
0& q+r > n \text{ or } r \leq 1\\
1& q+r=n \text{ and } n \geq q+2\\
g^0_{n-1,q,r} + \sum\limits_{i=r}^{n-q-1}g^0_{n-1,q,i} + \sum\limits_{i=r-1}^{n-q-1} g^1_{n-1,q,i}& \text{otherwise},
\end{cases}
\end{split}
\end{equation}

\begin{equation}
\begin{split}
&g^1_{n,q,r} =\\
&\quad \begin{cases}
0&q=0\\
0&q=n \text{ and } r=0\\
0&n=q=r=1\\
0&q+r>n\\
1&q+r=n, r>0\\
\sum_{i=q}^{q+r-1}g^1_{n-1,i,q+r-i}+ \sum_{i=q}^{q+r-2} g^0_{n-1,i,q+r-i}&q+r<n.
\end{cases}
\end{split}
\end{equation}

\subsection{The functional equations} \label{eq2}

It is clear from the recurrence that there is only one $(q-2,q,r)$ type node per level of the generating tree.  In particular, $g^2_{n,n+1,0}=1$ and $g^2_{n,q,r}=0$ if $q \neq n+1$ or $r \neq 0$.  Therefore $\frac{z}{1-z}$ is the generating function where the coefficient of $z^n$ is the number of $g^2_{n,q,r}$ type nodes summed over all values of $q$ and $r$.

For nodes of type $(q,q,r)$ and of type $(q-1,q,r)$, we first compute arrays for initial data.  Let $A0_n$ be the array where $g^0_{n,q,r}$ is in row $q$, column $r-1$, and let $A1_n$ be the array where $g^1_{n,q,r}$ appears in row $q$ column $r$.  We have the data shown in Table \ref{multitable}.

\begin{table}
\begin{center}
\begin{tabular}{|l|c|c|}
\hline
$n$&$A0_n$&$A1_n$\\
\hline
2&&$\left[\begin{array}{c}1\end{array}\right]$\\
\hline
3&$\left[\begin{array}{c}1\end{array}\right]$&$\left[\begin{array}{cc}1&1\\1&0\end{array}\right]$\\
\hline
4&$\left[\begin{array}{cc}4&1\\1&0\end{array}\right]$&$\left[\begin{array}{ccc}1&3&1\\1&1&0\\1&0&0\end{array}\right]$\\
\hline
5&$\left[\begin{array}{ccc}14&6&1\\4&1&0\\1&0&0\end{array}\right]$&$\left[\begin{array}{cccc}1&8&5&1\\1&3&1&0\\1&1&0&0\\1&0&0&0\end{array}\right]$\\
\hline
6&$\left[\begin{array}{cccc}50&27&8&1\\14&6&1&0\\4&1&0&0\\1&0&0&0\end{array}\right]$&$\left[\begin{array}{ccccc}1&23&19&7&1\\1&8&5&1&0\\1&3&1&0&0\\1&1&0&0&0\\1&0&0&0&0\end{array}\right]$\\
\hline
7&$\left[\begin{array}{ccccc}187&113&44&10&1\\50&27&8&1&0\\14&6&1&0&0\\4&1&0&0&0\\1&0&0&0&0\end{array}\right]$&$\left[\begin{array}{cccccc}1&74&69&34&9&1\\1&23&19&7&1&0\\1&8&5&1&0&0\\1&3&1&0&0&0\\1&1&0&0&0&0\\1&0&0&0&0&0\end{array}\right]$\\
\hline
8&$\left[\begin{array}{cccccc}730&468&212&65&12&1\\187&113&44&10&1&0\\50&27&8&1&0&0\\14&6&1&0&0&0\\4&1&0&0&0&0\\1&0&0&0&0&0\end{array}\right]$&$\left[\begin{array}{ccccccc}1&262&256&147&53&11&1\\1&74&69&34&9&1&0\\1&23&19&7&1&0&0\\1&8&5&1&0&0&0\\1&3&1&0&0&0&0\\1&1&0&0&0&0&0\\1&0&0&0&0&0&0\end{array}\right]$\\
\hline
\end{tabular}
\end{center}
\caption{Arrays containing $g_{n,q,q,r}$ and $g_{n,q-1,q,r}$ for small values of $n$}
\label{multitable}
\end{table}

Let $$C(x,y,z) = \sum_{n \geq 3} \sum_{q \geq 1} \sum_{r \geq 2} g^0_{n,q,r} x^q y^r z^n,$$

$$D(x,y,z) = \sum_{n \geq 2} \sum_{q \geq 1} \sum_{r \geq 1} g^1_{n,q,r} x^q y^r z^n.$$

From the recurrence for $g^0_{n,q,r}$, we have that

\begin{equation}
\begin{split}
C(x,y,z) = &zC(x,y,z)+\frac{zy}{y-1}C(x,y,z)-\frac{zy^2}{y-1}C(x,1,z)\\
& + \frac{xy^2z^2}{(1-z)(1-yz)(1-xz)}\\
&+\frac{zy^2}{y-1}\left(D(x,y,z) -\frac{xyz^2}{(1-xz)(1-yz)} \right)\\
&-\frac{zy^2}{y-1}\left(D(x,1,z)-\frac{xz^2}{(1-xz)(1-z)} \right).
\end{split}
\end{equation}

From the recurrence for $g^1_{n,q,r}$, we have that

\begin{equation}
\begin{split}
D(x,y,z) = &\frac{xyz^2}{(1-xz)(1-yz)}\\
& + \frac{zx}{x-y}D(x,y,z)-\frac{zx}{x-y}D(y,y,z)\\
&+\frac{zx}{x-y}C(x,y,z)-\frac{zx}{x-y}C(y,y,z).
\end{split}
\end{equation}

These equations are complicated to solve by hand, but the structure evident in $A0_n$ and $A1_n$ for small $n$ makes it easier to predict a solution by hand and verify that the solution works. In particular, notice that $A0_n$ looks like $A0_{n-1}$ with a column of 0s added to the right and a new row added to the top.  The same holds true for $A1_n$.  In other words the $(i,j)$ entry of $A0_n$ is the $(i-1,j)$ entry of $A0_{n-1}$ for $i \geq 2$.

It turns out that the entries of the first column of $A0_n$ are terms in $\frac{f(z)-1}{1-z}$ where $f(z) = \frac{1-z-\sqrt{1-6z+5z^2}}{2z} = 1+z+3z+10z^3+36z^4+137z^5+\cdots$ is the generating function for the first differences of the binomial convolution of the Catalan numbers.

Further, let $g(z) = \frac{16z^2(z-1)}{\left(1-z+\sqrt{1-6z+5z^2}\right)^3\left(-1+3z+\sqrt{1-6z+5z^2}\right)}$.  It also turns out that the generating function for column $i$ in $A0_n$ is the generating function for column $i-1$ times $g(z)$.

As for $A1_n$, we have the same effect of entry $(i,j)$ in $A1_n$ matching entry $(i-1,j)$ in $A1_{n-1}$, but the generating function is different.  It turns out that entry A091698 in the OEIS matches these entries but with additional minus signs.  From this structure, we may reverse engineer a conjecture for the entries in $A1_n$.

The unique solution to this system that matches the initial coefficients given above is

\begin{equation}
C(x,y,z)=\frac{2xy^2z^3}{(1-xz)((1-(y+1)z)\sqrt{5z^2-6z+1}+(1-(y+3)z)(1-z))},
\end{equation}

\begin{equation}
D(x,y,z)=\frac{2xyz^2}{(1-xz)(y\sqrt{5z^2-6z+1} +yz-2z-y+2)}.
\end{equation}

\subsection{The punchline}

We have determined closed forms for $C(x,y,z)$ and $D(x,y,z)$.  Since we really care only about the total number of nodes at level $n$, the total number of nodes where $p=q-2$ is 1 (generating function $\frac{z}{1-z}$), where $p=q-1$ is $D(1,1,z)$, and where $p=q$ is $C(1,1,z)$.  Adding and simplifying yields

\begin{equation}
C(1,1,z)+D(1,1,z)+\frac{z}{1-z} = \frac{-1+z+\sqrt{5z^2-6z+1}}{2(z-1)}
\end{equation}

This is, as per OEIS entry A007317, the generating function for the binomial convolution of the Catalan numbers.

\section{Conclusion}

In this paper we identified two sets, namely $\mathcal{A}_{201,210}(n)$ and $\mathcal{A}_{0021}(n)$, whose enumeration is given by sequence A007317 in OEIS.   Verifying this enumeration for 0021-avoiding ascent sequences completes the Wilf-classification of length 4 patterns for ascent sequences in conjunction with the work of Duncan and Steingr\'{i}msson \cite{DS11} and Mansour and Shattuck \cite{MS14}.  Both results use generating trees and experimentally-derived multivariate generating functions.  It remains open, and appears quite challenging, to find a statistic on each of these sets of ascent sequences corresponding to $k$ in the explicit formula of $\sum_{k=0}^{n-1}\binom{n-1}{k}C_k$ for the enumeration sequence.

\section*{Acknowledgement}
This paper grew out of a joint project with Andrew Baxter at Pennsylvania State University and the author is indebted to him for many helpful exposition comments.  Also, thank you to two anonymous referees who made a number of helpful presentation suggestions.

\begin{appendices}
\section*{Appendix }

\section{Proof of Lemma \ref{BigLemma}} \label{Appendix}

\begin{proof} We prove each part of Lemma \ref{BigLemma} in turn.

\medskip

\noindent\textbf{Part \ref{L7}:} Since $i \geq 2$, $i-1 \geq 1$. By Equation \ref{Recur1} we have  $d_{n,i} = g_{n,i-1,i} = g_{n-1,i-1,i} + \sum_{j=0}^{i-2}g_{n-1,j,i-1} = d_{n-1,i} + c_{n-1,i-1}$, as desired.

\medskip

\noindent\textbf{Part \ref{L5}:}  First, consider the case where $i=0$.  By Equation \ref{Recur1}, we have  $g_{n,0,n} = g_{n-1,0,n} + \sum_{i=1}^{n-q}i\cdot g_{n-1,i,n-1+i}$. In all terms $g_{n^*,p^*,q^*}$ on the right hand side, $q^*>n^*$, so all terms on the right are equal to 0. Therefore, $g_{n,0,n}=0$ for $n \geq 2$.

Now, if $i>0$, we proceed by induction.  For $n=1$, we have $g_{1,0,1}=1$, and $g_{1,i,1}=0$ for $i \neq 0$.  Next, assume that $g_{n-1,i,n-1}=0$ for $i \neq n-2$, and consider $g_{n,i,n}$ where $i \leq n-2$.  We have $g_{n,i,n} = g_{n-1,i,n} + \sum_{j=0}^{i-1}g_{n-1,j,n-1}$.  Since $n>n-1$, we have that $g_{n-1,i,n}=0$.  Also, since $j\leq i-1 \leq n-3$, by the induction hypothesis, we have that all $g_{n-1,j,n-1}$ terms in the sum are equal to 0.

\medskip

\noindent\textbf{Part \ref{L6}:}  We know that $g_{n,0,1}=1$, $d_{n,1} = g_{n,0,1}$, and $c_{n,1} = g_{n,0,1}$ by the definitions above.

We prove that $d_{n,n} = c_{n,n}=1$ by induction on $n$.  By part \ref{L5}, $c_{n,n} = \sum_{i=0}^{n-1} g_{n,i,n} = g_{n,n-1,n} = d_{n,n}$ since all terms except the last one in the sum are 0.

Now, for the base case, we see that $d_{1,1}=c_{1,1}=1$.

For the induction step, we assume that $d_{n-1,n-1}=c_{n-1,n-1}=1$, and we show that $d_{n,n}=c_{n,n}=1$.

We have: $d_{n,n} = g_{n,n-1,n} = g_{n-1,n-1,n} + \sum_{i=0}^{n-2}g_{n-1,i,n-1}$.  We know $g_{n-1,n-1,n}=0$ since $n>n-1$.  We know from part \ref{L5} that $\sum_{i=0}^{n-3}g_{n-1,i,n-1}=0$.  Therefore $d_{n,n} = g_{n-1,n-2,n-1} = d_{n-1,n-1}$, which is 1 by the induction hypothesis.

\medskip

Using parts \ref{L7} and \ref{L6}, we have completely characterized the $d_{n,i}$ terms in a recursive manner.  Parts \ref{L5} and \ref{L6} are used in the proof of parts \ref{L8} and \ref{L9}.  Notice also that parts \ref{L8} and \ref{L9} are symbiotic.  Part \ref{L8} determines the non-diagonal entries of $A_n$ that are not in the first row, while part \ref{L9} determines the non-diagonal entries of the first row.  Since parts \ref{L7}, \ref{L5}, and \ref{L6} determine the diagonal entries and the rightmost column of $A_n$, repeated applications of \ref{L8} and \ref{L9} completely determine the remaining entries of $A_n$.
 
\medskip

\noindent \textbf{Parts \ref{L8} and \ref{L9}:}  We prove these parts together by induction on $n$.  We assume that $A_{n}$ is characterized by the two parts for $n<n^*$ and use this assumption to show that the parts characterize the entries of $A_{n^*}$.

The statement of part \ref{L8} requires $n \geq 4$.  In this case, $p=1$ and $q=3$ are the only values that satisfy the inequality.  Using parts \ref{L7}, \ref{L5}, and \ref{L6} and Equation \ref{Recur1}, we see 
$g_{4,1,3} = g_{3,1,3} + g_{3,0,2} = 0 + g_{3,0,2} = g_{2,0,2} + g_{2,1,2} = 0 + 1 = 1$ and $g_{4,3,4}=1$. 

In part \ref{L9}, we need $n \geq 3$.  It is quickly verified that $g_{3,0,2} = 1$ and $g_{3,0,3} + g_{3,1,3} + g_{3,2,3} = 0+0+1=1$. Also for $n=4$, we have $g_{4, 0, 2} = 5$ and $g_{4,0,3}+g_{4,1,3}+g_{4,2,3}=0+1+4=5$.  Also, $g_{4,0,3}=0$ and $g_{4,0,4}+g_{4,1,4}+g_{4,2,4}=0+0+0=0$.  Therefore, part \ref{L8} and \ref{L9} hold for $n \leq 4$.

Next, we assume that both parts simultaneously hold for $n < n^*$ and consider the entries of $A_{n^*}$.

For part \ref{L8}, and consider $g_{n^*,p,q}$ where $3 \leq p+2 < q+1 \leq n^*$.  Notice that $3 < q+1 \leq n^*$ indicates $2 < q < n^*-1$ so part \ref{L9} also applies. From Equation \ref{Recur1}, $$g_{n^*,p,q} = g_{n^*-1,p,q}+ \sum_{i=0}^{p-1} g_{n^*-1,i,q-1}.$$  By the induction hypothesis of part \ref{L8} applied to every term except for $g_{n^*-1,0,q-1}$ and the induction hypothesis of part \ref{L9} applied to this remaining term, we have: 

\begin{equation}
\begin{split}
g_{n^*,p,q} &= g_{n^*-1,p,q}+ \sum_{i=0}^{p-1} g_{n^*-1,i,q-1} \\
& = g_{n^*-1,p+2,q+1} + g_{n^*-1,0,q-1} + \sum_{i=1}^{p-1} g_{n^*-1,i+2,q}\\
& = g_{n^*-1,p+2,q+1} + g_{n^*-1,0,q-1} + \sum_{i=3}^{p+1} g_{n^*-1,i,q}\\
&= g_{n^*-1,p+2,q+1} + \left(g_{n^*-1,0,q} +g_{n^*-1,1,q}+g_{n^*-1,2,q}\right)+\sum_{i=3}^{p+1} g_{n^*-1,i,q}\\
&= g_{n^*-1,p+2,q+1} +\sum_{i=0}^{p+1} g_{n^*-1,i,q}\\
& = g_{n^*,p+2,q+1}.\\
\end{split}
\end{equation}

For part \ref{L9}, consider $g_{n^*,0,q+1} + g_{n^*,1,q+1} + g_{n^*,2,q+1}$ where $2 \leq q \leq n^*-1$.

By Equation \ref{Recur1} we have:

\begin{equation}
\begin{split}
g_{n^*,0,q+1}& + g_{n^*,1,q+1} + g_{n^*,2,q+1}=\\
&\left(g_{n^*-1,0,q+1} + \sum_{i=1}^{n^*-q-1}i \cdot g_{n^*-1,i,q+i} \right)\\
&\quad + \left(g_{n^*-1,1,q+1} + g_{n^*-1,0,q}\right)\\
&\quad + \left(g_{n^*-1,2,q+1} + g_{n^*-1,0,q} + g_{n^*-1,1,q}\right).\\
\end{split}
\end{equation}

We wish to show that this quantity is equal to $g_{n^*,0,q} = g_{n^*-1,0,q} + \sum_{i=1}^{n^*-q} i\cdot g_{n^*-1,i,q-1+i}$.

First, apply part \ref{L8} to the summation to obtain
\begin{equation}
\begin{split}
g_{n^*,0,q+1}& + g_{n^*,1,q+1} + g_{n^*,2,q+1}=\\
& \left(g_{n^*-1,0,q+1} + \sum_{i=1}^{n^*-q-1}i \cdot g_{n^*-1,i+2,q+1+i} \right)\\
&\quad + \left(g_{n^*-1,1,q+1} + g_{n^*-1,0,q}\right)\\
&\quad + \left(g_{n^*-1,2,q+1} + g_{n^*-1,0,q} + g_{n^*-1,1,q}\right).\\
\end{split}
\end{equation}

Next, notice that we may apply part \ref{L9} repeatedly to a term of the form $g_{a,0,b}$ to obtain
\begin{equation}
\begin{split}
g_{a,0,b} &= g_{a,0,b+1}+g_{a,1,b+1}+g_{a,2,b+1}\\
&=\left(g_{a,0,b+2}+g_{a,1,b+2}+g_{a,2,b+2}\right)+g_{a,1,b+1}+g_{a,2,b+1}\\
&=\left(\left(g_{a,0,b+3}+g_{a,1,b+3}+g_{a,2,b+3}\right)+g_{a,1,b+2}+g_{a,2,b+2}\right)+g_{a,1,b+1}+g_{a,2,b+1}\\
&=\cdots\\
&=g_{a,0,a} + \sum_{i=1}^{a-b} \left(g_{a,1,b+i} + g_{a,2,b+i}\right).\\
\end{split}
\end{equation}

Apply part \ref{L9} repeatedly to both the $g_{n^*-1,0,q+1}$ term and one of the $g_{n^*-1,0,q}$ terms, and combine like terms to obtain
\begin{equation}
\begin{split}
g_{n^*,0,q+1}& + g_{n^*,1,q+1} + g_{n^*,2,q+1}\\
&= g_{n^*-1,0,q} + g_{n^*-1,1,q} + \sum_{i=1}^{n^*-q-1}i \cdot g_{n^*-1,i+2,q+1+i}\\
& + \sum_{i=1}^{n^*-1}2g_{n^*-1,1,q+i} + \sum_{i=1}^{n^*-1}2g_{n^*-1,2,q+i}.\\
\end{split}
\end{equation}

Now, given $g_{n^*-1,a,q+i}$ with $a \in \{1,2\}$, we know that $$g_{n^*-1,a,q+i} = g_{n^*-1,a+2j,q+i+j}$$ for all $j \geq 1$.  Further, given, $a$, $i$, and $q$, there is a unique value $j^*$ such that $(q+i+j^*)-(a+2j^*) = q-1$.  Indeed, $j^*=i+1-a$.  Rewrite each term $g_{n^*-1,a,q+i}$ as $g_{n^*-1,a+2j^*,q+i+j^*}$.

Further, notice that given $g_{n^*-1,a+2j^*,q+i+j^*}$, there is a unique pair $(a,i)$ with $a \in \{1,2\}$ and $i \geq 1$ that produces $(a+2j^*,q+i+j^*)$.  Since $a \in \{1,2\}$, the value of $a$ is determined by the parity of $a+2j^*$.  This determines the value of $j^*$, and since $q$ is fixed, the value of $j^*$ determines the value of $i$.  Therefore, 
\begin{equation}
\begin{split}\sum_{i=1}^{n^*-1}2g_{n^*-1,1,q+i} + \sum_{i=1}^{n^*-1}2g_{n^*-1,2,q+i} = \sum_{i=2}^{n^*-q} 2g_{n^*-1 , i, q-1+i}.
\end{split}
\end{equation}
  We have:

\begin{equation}
\begin{split}
g_{n^*,0,q+1}& + g_{n^*,1,q+1} + g_{n^*,2,q+1}\\
&= g_{n^*-1,0,q} + g_{n^*-1,1,q} + \sum_{i=1}^{n^*-q-1}i \cdot g_{n^*-1,i+2,q+1+i} +\sum_{i=2}^{n^*-q} 2g_{n^*-1 , i, q-1+i}\\
&= g_{n^*-1,0,q} + g_{n^*-1,1,q} + \sum_{i=2}^{n^*-q}(i-2) \cdot g_{n^*-1,i,q-1+i} +\sum_{i=2}^{n^*-q} 2g_{n^*-1 , i, q-1+i}\\
&= g_{n^*-1,0,q} + g_{n^*-1,1,q} + \sum_{i=2}^{n^*-q}(i) \cdot g_{n^*-1,i,q-1+i}\\
&= g_{n^*-1,0,q} + \sum_{i=1}^{n^*-q}(i) \cdot g_{n^*-1,i,q-1+i}\\
&=g_{n^*,0,q}.\\
\end{split}
\end{equation}

\medskip

The final two parts of Lemma \ref{BigLemma} are more straightforward.  Part \ref{L10} is a statement about partial column sums of $A_n$ that follows directly from parts \ref{L8} and \ref{L9}.  

\medskip

\noindent\textbf{Part \ref{L10}:} We have:

\begin{equation}
\begin{split}
\sum_{i=0}^k g_{n,i,q}  &= g_{n,0,q} + \sum_{i=1}^k g_{n,i,q}\\
& = g_{n,0,q+1}+g_{n,1,q+1} + g_{n,2,q+1} + \sum_{i=1}^k g_{n,i+2,q+1}\\
& = \sum_{i=0}^{k+2} g_{n,i,q+1}.\\
\end{split}
\end{equation}

\medskip

\noindent\textbf{Part \ref{L11}:} This is a direct consequence of part \ref{L10} that provides another linear relationship between the $d_{n,i}$ and $c_{n,i}$ terms.  Take part \ref{L10} with $k=q-2$. Then
$$\sum_{j=0}^{q-2} g_{n,j,q} = \sum_{j=0}^{q} g_{n,j,q+1}.$$

The right hand side is $c_{n,q+1}$, while the left is $c_{n,q} - d_{n,q}$.  Let $q=i-1$ to see the statement holds.

\end{proof}

\end{appendices}

\end{document}